\documentclass{amsart}
\usepackage{amssymb,color,hhline}
\usepackage{amsmath}
\usepackage{a4}
\def\nmonth{\ifcase\month\ \or Enero\or
Febrero\or Marzo\or Abril\or Mayo\or Junio\or Julio\or Agosto\or
Septiembre\or Octubre\or Noviembre\else Diciembre\fi}
\newtheorem{theorem}{Theorem}[section]
\newtheorem{lemma}[theorem]{Lemma}

\makeatletter
  
  \@addtoreset{equation}{section}
 \makeatother

\newcommand{\g}[2]{\ensuremath{\langle #1,#2 \rangle}}
\newcommand{\ja}{\ensuremath{\mathcal{J}}}

\begin{document}

\title[The complex Jacobi operator]{Relating the curvature tensor and the complex Jacobi operator of an almost
Hermitian manifold}
\author{M. Brozos-V{\'a}zquez, E. Garc{\'\i}a-R{\'\i}o and P. Gilkey}
\begin{address}{M. Brozos-V{\'a}zquez, Department of Geometry and Topology, Faculty of
Mathematics, University of Santiago de Compostela, Santiago de Compostela
15782, Spain.}\end{address}
\begin{email}{mbrozos@usc.es}\end{email}
\begin{address}{E. Garc{\'\i}a-R{\'\i}o, Department of Geometry and Topology, Faculty of
Mathematics, University of Santiago de Compostela, Santiago de Compostela
15782, Spain.}\end{address}
\begin{email}{xtedugr@usc.es}\end{email}
\begin{address}{P. Gilkey, Mathematics Department, University of Oregon, Eugene Or
97403, USA}
\end{address}\begin{email}{gilkey@uoregon.edu}\end{email}

\begin{abstract} Let $J$ be a unitary almost complex structure on a
Riemannian manifold $(M,g)$. If $x$ is a unit tangent vector, let
$\pi:=\operatorname{Span}\{x,Jx\}$ be the associated complex line
in the tangent bundle of $M$. The complex Jacobi operator and the
complex curvature operators are defined, respectively, by
$\ja(\pi):=\ja(x)+\ja(Jx)$ and
$\mathcal{R}(\pi):=\mathcal{R}(x,Jx)$. We show that if $(M,g)$ is
Hermitian or if $(M,g)$ is nearly K\"ahler, then either the
complex Jacobi operator or the complex curvature operator
completely determine the full curvature operator; this generalizes
a well known result in the real setting to the complex setting. We
also show this result fails for general almost Hermitian
manifolds.\end{abstract}
\maketitle

\section{Introduction}
We shall let $\mathcal{M}:=(M,g)$ be a Riemannian manifold of
dimension $m$. Let
\begin{eqnarray*}
&&\mathcal{R}(x,y):=\nabla_x\nabla_y-\nabla_y\nabla_x-\nabla_{[x,y]},\\
&&R(x,y,z,w):=g(\mathcal{R}(x,y)z,w)
\end{eqnarray*}
be the curvature operator and the curvature tensor, respectively;
$R$ has the symmetries:
\begin{eqnarray}
&&R(x,y,z,w)=R(z,w,x,y)=-R(y,x,z,w),\label{symmetries-curv-tensor}\\
&&R(x,y,z,w)+R(y,z,x,w)+R(z,x,y,w)=0\,.\nonumber
\end{eqnarray}

It is convenient to work in the algebraic setting. Let $V$ be a
vector space of dimension $m$. We say that $A\in\otimes^4(V^*)$ is
an {\it algebraic curvature tensor} if $A$ has the symmetries
given in Equation (\ref{symmetries-curv-tensor}). We consider a
{\it model} $\mathfrak{M}:=(V,\langle\cdot,\cdot\rangle,A)$ where
$\langle\cdot,\cdot\rangle$ is an auxiliary positive definite
inner product on $V$. Every model is geometrically realizable;
given a model $\mathfrak{M}$, one can construct a Riemannian
manifold $\mathcal{M}$ so that $\mathfrak{M}$ is isomorphic to
$(T_PM,g_P,R_P)$ for some point $P\in M$.

Given a model $\mathfrak{M}$, one uses the inner product
$\langle\cdot,\cdot\rangle$ to raise indices and define an
associated curvature operator $\mathcal{A}$. The {\it Jacobi
operator} $\mathcal{J}:y\rightarrow\mathcal{A}(y,x)x$ is
characterized  by the identity
$$\langle\mathcal{J}(x)y,z\rangle=A(y,x,x,z)\,.$$

The Jacobi operator determines the full curvature tensor. The
following theorem is well known; Assertion (2) in the geometric
setting is an immediate consequence of the corresponding Assertion
(1) in the algebraic setting:

\begin{theorem}\label{thm-1.1}
\ \begin{enumerate} \item Let
$\mathfrak{M}_i=(V_i,\langle\cdot,\cdot\rangle_i,A_i)$ be models
for $i=1,2$. Suppose there exists an isometry
$\theta:(V_1,\langle\cdot,\cdot\rangle_1)\rightarrow
 (V_2,\langle\cdot,\cdot\rangle_2)$
so that $\mathcal{J}_{\mathfrak{M}_2}(\theta
x)\theta=\theta\mathcal{J}_{\mathfrak{M}_1}(x)$ for all $x\in
V_1$. Then $\theta^*A_2=A_1$. \item Let $\mathcal{M}_i=(M_i,g_i)$
be Riemannian manifolds for $i=1,2$. Suppose there is an isometry
$\theta:(T_PM_1,g_1)\rightarrow(T_QM_2,g_2)$ so that
$\mathcal{J}_{\mathcal{M}_2,Q}(\theta
x)\theta=\theta\mathcal{J}_{\mathcal{M}_1,P}(x)$ for all $x\in
T_PM_1$. Then $\theta^*R_2(Q)=R_1(P)$.
\end{enumerate}
\end{theorem}

If $J$ is a unitary almost complex structure on a Riemannian
manifold $(M,g)$, then $\mathcal{C}:=(M,g,J)$ is said to be an
{\it almost Hermitian manifold}. In the algebraic setting,
$\mathfrak{C}:=(V,\langle\cdot,\cdot\rangle,J,A)$ is said to be a
{\it complex model} if $J$ is a unitary complex structure and if
$A$ is an algebraic curvature tensor. Any point $P$ of an almost
Hermitian manifold $\mathcal{C}$ determines a corresponding
complex model $\mathfrak{C}(\mathcal{C},P):=(T_PM,g_P,J_P,R_P)$ in
a natural fashion.

Let $\mathfrak{C}$ be a complex model. The Ricci tensor $\rho$ and
the $\star$-Ricci tensor $\rho^\star$ are defined by contracting
indices. If $\{e_1,...,e_m\}$ is an orthonormal basis for $V$,
then
$$
\rho(x,y)=\sum_{i=1}^m R(e_i,x,y,e_i)\quad\text{and}\quad
\rho^\star (x,y):=\sum_{i=1}^m R(x,Je_i,Jy,e_i)\,.
$$
We note that $\rho^\star$ is not in general a symmetric
$2$-tensor; however one does have that
$\rho^\star(x,y)=\rho^\star(y,x)$ if the compatibility condition
given below in Lemma \ref{lemma-compatibility} is satisfied. The
scalar curvature $\tau$ and the $\star$-scalar curvature
$\tau^\star$ are defined by a final contraction:
$$\tau=\sum_{i=1}^m \rho(e_i,e_i)\quad\text{and}\quad
\tau^\star=\sum_{i=1}^m\rho^\star(e_i,e_i)\,.$$

 We say a
$2$-dimensional subspace $\pi$ of $V$ is a {\it complex line} if
$J\pi=\pi$. Let $\mathbb{CP}(V,J)$ be the complex projective space
of complex lines in $V$. If $\pi\in\mathbb{CP}(V,J)$, let $S(\pi)$
be the set of unit vectors in $\pi$. Let $x\in S(\pi)$ for
$\pi\in\mathbb{CP}(V,J)$. Then one has
$\pi=\pi_x:=\operatorname{Span}\{x,Jx\}$.  The {\it holomorphic
sectional curvature} $Q(\pi)$, {\it complex Jacobi operator}
$\mathcal{J}(\pi)$, and {\it complex skew-symmetric curvature
operator} $\mathcal{R}(\pi)$ are then defined for
$\pi\in\mathbb{CP}(V,J)$, respectively, by setting
\begin{equation}\label{hol-sect-curv}
\begin{array}{ll}
Q(\pi):=A(x,Jx,Jx,x),&\mathcal{J}(\pi):=\mathcal{J}(x)+\mathcal{J}(Jx),\\
\mathcal{R}(\pi):=\mathcal{R}(x,Jx)&\text{for any }x\in
S(\pi)\,.\vphantom{\vrule height 11pt}
\end{array}\end{equation}

We shall be considering several important families of almost
Hermitian manifolds. $\mathcal{C}$ is said to be {\it Hermitian}
if the Nijenhuis tensor vanishes, i.e. if
\[
[X,Y]+J[JX,Y]+J[X,JY]-[JX,JY]=0\quad\text{for all}\quad X,Y\,.
\]
 Equivalently, see \cite{newlander-niremberg}, this
means that we can find local holomorphic coordinates
$z_\nu=x_\nu+\sqrt{-1}y_\nu$ for $1\le\nu\le\frac12m$ so that
$J\partial_{x_\nu}=\partial_{y_\nu}$ and
$J\partial_{y_\nu}=-\partial_{x_\nu}$; the transition functions
relating two such coordinate systems are then complex analytic.
There are other natural assumptions that may be imposed and which
define other important families. For example, one says that
$\mathcal{C}$ is {\it nearly K{\"a}hler} if $(\nabla_xJ)x=0$ for all
tangent vectors $x$; we refer to \cite{nagy} for further
information concerning this class of manifolds. We say that
${\mathcal{C}}$ is {\it almost K{\"a}hler} if the two form
$\Omega(x,y):=\g{Jx}{y}$ is closed; we refer to
\cite{apostolov-draghici} for a survey and to
\cite{apostolov-armstrong-draghici}, \cite{kir} and
\cite{oguro-sekigawa} for some recent results concerning this
class of manifolds.

The following result, which generalizes Theorem \ref{thm-1.1} to
the complex setting, is the central result of this paper. It shows
that the full curvature tensor is determined either by the complex
Jacobi operator or by the complex curvature operator in certain
natural geometric contexts:

\begin{theorem}\label{complex-dif-theorem}
For $i=1,2$, let $\mathcal{C}_i=(M_i,g_i,J_i)$ be either Hermitian or nearly
K{\"a}hler manifolds. Let
$\theta:(T_PM_1,g_1,J_1)\longrightarrow(T_QM_2,g_2,J_2)$ be a
complex isometry. The following assertions are equivalent:
\begin{enumerate}
\item $\theta
\mathcal{J}_{\mathcal{C}_1,P}(\pi)=\mathcal{J}_{\mathcal{C}_2,Q}(\theta\pi)\theta$
for all $\pi\in\mathbb{CP}(T_PM_1,J_{1,P})$. \item $\theta
R_{\mathcal{C}_1,P}(\pi)=R_{\mathcal{C}_2,Q}(\theta\pi)\theta$ for
all $\pi\in\mathbb{CP}(T_PM_1,J_{1,P})$. \item $\theta^\ast
R_{\mathcal{C}_2,Q}=R_{\mathcal{C}_1,P}$\,.
\end{enumerate}
\end{theorem}

Our result for almost K\"ahler setting manifolds is a bit weaker:

\begin{theorem}\label{thm-1.3}
Let $\mathcal{C}=(M,g,J)$ be an almost K{\"a}hler manifold. Assume
either that $\ja(\pi)=0$ for all $\pi\in\mathbb{CP}(TM,J)$ or that
$\mathcal{R}(\pi)=0$ for all $\pi\in\mathbb{CP}(TM,J)$. Then $M$
is flat.
\end{theorem}

Let $\mathfrak{C}$ be a complex model. There is a basic
compatibility condition we work with that relates the structures
$J$ and $A$:

\begin{lemma}\label{lemma-compatibility}
Let $\mathfrak{C}=(V,\langle\cdot,\cdot\rangle,J,A)$ be a complex
model. The following conditions are equivalent and if any is
satisfied, we shall say that $\mathfrak{C}$ is a {\rm compatible
complex model}.
\begin{enumerate}
\item $J^*A=A$, i.e. $A(x,y,z,t)=A(Jx,Jy,Jz,Jt)$ for all
$x,y,z,t\in V$. \item $\mathcal{J}(\pi)J=J\mathcal{J}(\pi)$ for
all $\pi\in\mathbb{CP}(V,J)$. \item
$\mathcal{A}(\pi)J=J\mathcal{A}(\pi)$ for all
$\pi\in\mathbb{CP}(V,J)$.
\end{enumerate}\end{lemma}

We note, see Lemma \ref{lem-3.1}, that if $\mathcal{C}=(M,g,J)$ is
a nearly K{\"a}hler manifold, then
$\mathfrak{C}(\mathcal{C},P)=(T_PM,g_P,J_P,R_P)$ is a compatible
complex model for any point $P\in M$. In general, a manifold
satisfying this condition at every point is known in the
literature as a $RK$-manifold.

What is perhaps rather surprising is that Theorem
\ref{complex-dif-theorem} does not have a purely algebraic
analogue even if we impose the compatibility condition of Lemma
\ref{lemma-compatibility}:
\begin{theorem}\label{thm-1.5}
If $m\equiv0$ mod $4$, there exists a complex model
$\mathfrak{C}=(V,\langle\cdot,\cdot\rangle,J,A)$ with $A\ne 0$ so
that
 $\mathcal{J}(\pi)=0$ and so that
$\mathcal{A}(\pi)=0$ for every $\pi\in\mathbb{CP}(V,J)$.
\end{theorem}

Let $\mathcal{C}=(M,g,J)$ be an almost Hermitian manifold. We let
$\mathcal{U}_{\mathcal{C}}$ be the bundle of complex isometries of
$TM$; the fibers of this bundle are the associated unitary group
of the fibers. If $\Theta\in
C^\infty\{\mathcal{U}_{\mathcal{C}}\}$ and if $P\in M$, then
$\theta_P:=\Theta(P)$ is a complex isometry of $(T_PM,g_P,J_P)$
for any $P\in M$. We show that Theorem \ref{complex-dif-theorem}
fails in the almost Hermitian context by establishing the
following result:

\begin{theorem}\label{thm-1.6}
Let $m\equiv0$ mod $4$. There exists an almost Hermitian manifold
$\mathcal{C}$ and there exists $\Theta\in
C^\infty\{\mathcal{U}_{\mathcal{C}}\}$ so that for any point $P$
in $M$ we have:
\begin{enumerate}
\item
$\theta_P\mathcal{J}_{\mathcal{C}}(\pi)=\mathcal{J}_{\mathcal{C}}(\theta_P\pi)\theta_P$
 for all $\pi\in\mathbb{CP}(T_PM,J)$.
\item
$\theta_P\mathcal{R}_{\mathcal{C}}(\pi)=\mathcal{R}_{\mathcal{C}}(\theta_P\pi)\theta_P$
 for all $\pi\in\mathbb{CP}(T_PM,J)$.
\item $\theta_P^*R_P\ne R_P$.
\end{enumerate}\end{theorem}

Here is a brief outline to this paper. Section \ref{sect-alg} is
algebraic in nature. We begin by establishing Lemma
\ref{lemma-compatibility}. Next, in Lemma \ref{vanhecke-lemma}, we
prove a result of Vanhecke \cite{vanhecke} which expresses
$R(x,y,y,x)$ for a compatible complex model in terms of the
holomorphic sectional curvature $Q$ defined in Equation
(\ref{hol-sect-curv}) and in terms of an additional tensor
\begin{equation}\label{lambda}
\lambda(x,y)=R(x,y,y,x)-R(x,y,Jy,Jx)\,.
\end{equation}
The identity of Lemma \ref{vanhecke-lemma} is polarized to
establish a result of Sato \cite{sato} in Lemma \ref{sato-lemma}.
We then turn to a study of the complex Jacobi operator by studying
the condition $\mathcal{J}(\pi)=0$ for all
$\pi\in\mathbb{CP}(V,J)$ in Lemma \ref{lem-2.3} and show this
implies that $\rho$ and $\rho^*$ both vanish. We then prove
Theorem \ref{thm-1.5}. We conclude Section \ref{sect-alg} by
relating Lemma \ref{lem-2.3} to the curvature decompositions of
Gray \cite{gray} in Lemma \ref{lem-2.4}.

In Section \ref{sect-geo} we apply the results of Section
\ref{sect-alg} to the geometric context. We begin by recalling
certain results of Gray and Yano concerning Hermitian, nearly
K{\"a}hler and almost K{\"a}hler manifolds. These results are then applied
to prove Theorems \ref{complex-dif-theorem} and \ref{thm-1.3}. The
construction used to establish Theorem \ref{thm-1.5} in the
algebraic setting is then used to prove Theorem \ref{thm-1.6} in
the geometric setting.

There are many example in the literature of almost K{\"a}hler
manifolds which are not K{\"a}hler
\cite{abbena,cordero-marisa-mdeleon,watson}. Also, a great
interest has been shown in finding conditions for an almost K{\"a}hler
manifold to be K{\"a}hler. For example, the Goldberg conjecture
states: {\it A compact Einstein almost K{\"a}hler manifold is K{\"a}hler}.
This conjecture has generated extensive literature, see for
example \cite{sekigawa-1987,sekigawa-vanhecke}. Many other
conditions have been studied which might imply that an almost
K{\"a}hler manifold is K{\"a}hler, see
\cite{apostolov-armstrong-draghici,balas-gauduchon} for example.
We conclude the paper in Section \ref{sect-add-res} with two
related results.

We shall adopt the following notational conventions. The curvature
tensor and curvature operator of a Riemannian manifold will be
denoted by $R$ and $\mathcal{R}$, respectively; an algebraic
curvature tensor and the corresponding algebraic curvature
operator will be denoted by $A$ and $\mathcal{A}$, respectively. A
real Riemannian manifold and a real model will be denoted by
$\mathcal{M}=(M,g)$ and
$\mathfrak{M}=(V,\langle\cdot,\cdot\rangle,A)$, respectively. An
almost Hermitian manifold and a complex model will be denoted by
$\mathcal{C}=(M,g,J)$ and
$\mathfrak{C}=(V,\langle\cdot,\cdot\rangle,J,A)$, respectively.
The Jacobi operator will be denoted by $\mathcal{J}$; we will
subscript as appropriate when more than one curvature operator is
under consideration.

In this preprint, we have chosen to give full details of many algebraic
computations in the interests of keeping matters as self-contained as possible
for the convenience of the reader; the version that will be submitted for
publication will be a bit shorter as we shall omit many of these computations
if they are available elsewhere.

\section{Algebraic Results}\label{sect-alg}

\begin{proof} We first establish Lemma \ref{lemma-compatibility}.
Suppose first that Assertion (1) holds, i.e. that
$A(x,y,z,t)=A(Jx,Jy,Jz,Jt)$ for all $x$, $y$, $z$, $t$. Then
\begin{equation}\label{eqn-2.b}
A(y,x,x,Jz)+A(y,Jx,Jx,Jz)=-A(Jy,x,x,z)-A(Jy,Jx,Jx,z)\,,
\end{equation}
which implies $\g{\ja(\pi_x)y}{Jz}=-\g{\ja(\pi_x)Jy}{z}$ and,
hence, $\ja(\pi_x)J=J\ja(\pi_x)$. Assume conversely that
$\ja(\pi_x)J=J\ja(\pi_x)$ or equivalently that Equation
(\ref{eqn-2.b}) holds for all $x$. Polarizing this identity and
replacing $z$ by $-Jz$ yields
\begin{eqnarray}
&&A(y,x,w,z)+A(y,w,x,z)+A(y,Jx,Jw,z)+A(y,Jw,Jx,z)\nonumber\\
&=&A(Jy,x,w,Jz)+A(Jy,w,x,Jz)+A(Jy,Jx,Jw,Jz)\label{eqn-2.c}\\
&+&A(Jy,Jw,Jx,Jz)\,.\nonumber
\end{eqnarray}
Interchanging arguments $1 \leftrightarrow 2 $ and $3
\leftrightarrow 4 $ in the curvature tensors then yields:
\begin{eqnarray*}
&&A(x,y,z,w)+A(w,y,z,x)+A(Jx,y,z,Jw)+A(Jw,y,z,Jx)\\
&=&A(x,Jy,Jz,w)+A(w,Jy,Jz,x)+A(Jx,Jy,Jz,Jw)\\&+&A(Jw,Jy,Jz,Jx)\,.
\end{eqnarray*}

If we interchange $x$ and $y$ and we interchange $z$ and $w$ in
this identity, we get
\begin{eqnarray}
&&A(y,x,w,z)+A(z,x,w,y)+A(Jy,x,w,Jz)+A(Jz,x,w,Jy)\nonumber\\
&=&A(y,Jx,Jw,z)+A(z,Jx,Jw,y)+A(Jy,Jx,Jw,Jz)\label{eqn-2.d}\\
&+&A(Jz,Jx,Jw,Jy)\,.\nonumber
\end{eqnarray}

Adding (\ref{eqn-2.c}) and (\ref{eqn-2.d}) and simplifying yields:
\begin{equation}\label{eqn-2.e}
\phantom{=}A(y,x,w,z)+A(y,w,x,z) =A(Jy,Jx,Jw,Jz)+A(Jy,Jw,Jx,Jz)
\end{equation}
We permute the indices in Equation (\ref{eqn-2.e}) to change
$y\rightarrow x\rightarrow w\rightarrow y$. This yields:
\begin{equation}\label{eqn-2.f}
\phantom{=}A(x,w,y,z)+A(x,y,w,z) =A(Jx,Jw,Jy,Jz)+A(Jx,Jy,Jw,Jz)\,.
\end{equation}
We add 2(\ref{eqn-2.e}) and (\ref{eqn-2.f}) and use the Bianchi
identity to see
\begin{eqnarray*}
&&3A(y,w,x,z)=A(y,x,w,z)+2A(y,w,x,z)+A(x,w,y,z)\\
&=&A(Jy,Jx,Jw,Jz)+2A(Jy,Jw,Jx,Jz)+A(Jx,Jw,Jy,Jz)\\
&=&3A(Jy,Jw,Jx,Jz)\,.
\end{eqnarray*}
The desired identity now follows.

We now prove that Assertion (1) implies Assertion (3). We have:
\begin{eqnarray*}
&&\langle J\mathcal{A}(\pi_x)y,z\rangle=-\langle\mathcal{A}(\pi_x)y,Jz\rangle=-A(x,Jx,y,Jz)\\
&=&-A(Jx,JJx,Jy,JJz)=A(x,Jx,Jy,z)=\langle\mathcal{A}(\pi_x)Jy,z\rangle\,.
\end{eqnarray*}
Thus $J\mathcal{A}(\pi_x)=\mathcal{A}(\pi_x)J$ as desired.

We finally show that Assertion (3) implies Assertion (1). We have
\begin{eqnarray*}
&&J\mathcal{A}(x,Jx)=\mathcal{A}(x,Jx)J,\\
&\Rightarrow&\langle JA(x,Jx)z,w\rangle-\langle A(x,Jx)Jz,w\rangle=0,\\
&\Rightarrow&A(x,Jx,z,Jw)+A(x,Jx,Jz,w)=0\,.
\end{eqnarray*}
Polarizing yields an identity for all $x,y,z,w$:
\begin{eqnarray*}
0&=&A(y,Jx,z,Jw)+A(x,Jy,z,Jw)+A(y,Jx,Jz,w)\\
&+&A(x,Jy,Jz,w)\,.
\end{eqnarray*}
Interchange the first two arguments in the first and third term to
see:
\begin{eqnarray*}
0&=&-A(Jx,y,z,Jw)+A(x,Jy,z,Jw)-A(Jx,y,Jz,w)\\&+&A(x,Jy,Jz,w)\,.
\end{eqnarray*}
Replace $(x,w)$ by $(Jx,Jw)$ to show:
\begin{equation}\label{eqn-2.g}
\begin{array}{rcl}
0&=&-A(x,y,z,w)-A(Jx,Jy,z,w)+A(x,y,Jz,Jw)\\
&+&A(Jx,Jy,Jz,Jw)\,.
\end{array}\end{equation}
Interchange the first two arguments with the final two arguments:
\begin{eqnarray*}
0&=&-A(z,w,x,y)-A(z,w,Jx,Jy)+A(Jz,Jw,x,y)\\
&+&A(Jz,Jw,Jx,Jy)\,.
\end{eqnarray*}
Change notation to interchange $x$ and $z$, and $y$ and $w$, to
see:
\begin{equation}\label{eqn-2.h}
\begin{array}{rcl}
0&=&-A(x,y,z,w)-A(x,y,Jz,Jw)+A(Jx,Jy,z,w)\\
&+&A(Jx,Jy,Jz,Jw)\,.
\end{array}\end{equation}
\smallbreak\noindent We add Equations (\ref{eqn-2.g}) and
(\ref{eqn-2.h}) to conclude
$$-A(x,y,z,w)+A(Jx,Jy,Jz,Jw)=0$$
and complete the proof that Assertion (3) implies Assertion (1).
\end{proof}

The following result is due to Vanhecke \cite{vanhecke}; it was
originally stated in a purely geometrical setting. Let $Q$ and
$\lambda$ be defined by Equations (\ref{hol-sect-curv}) and
(\ref{lambda}), respectively.
\begin{lemma}\label{vanhecke-lemma}
Let $\mathfrak{C}=(V,\g\cdot\cdot,A,J)$ be a compatible complex
model. Then
\begin{eqnarray*}
32A(x,y,y,x)&=&3Q(x+Jy)+3Q(x-Jy)-Q(x+y)
-Q(x-y)\\&-&4Q(x)-4Q(y)+4\{5\lambda(x,y)+\lambda(x,Jy)\}\,.
\end{eqnarray*}
\end{lemma}
\begin{proof} First note that
\[
\begin{array}{rcl}
Q(x+y)&=&A(x+y,Jx+Jy,Jx+Jy,x+y)\\
\noalign{\smallskip}
&=&A(x,Jx,Jx,x)+A(x,Jy,Jy,x)\\
\noalign{\smallskip}&&+A(y,Jx,Jx,y)+A(y,Jy,Jy,y)\\
\noalign{\smallskip}
&&+2A(x,Jx,Jx,y)+2A(x,Jy,Jy,y)+2A(x,Jx,Jy,x)\\
\noalign{\smallskip} &&+2A(y,Jx,Jy,y)+2A(x,Jx,Jy,y)+2A(x,Jy,Jx,y).
\end{array}
\]
Hence \medbreak\quad$
3Q(x+Jy)+3Q(x-Jy)-Q(x+y)-Q(x-y)-4Q(x)-4Q(y)$\smallbreak\qquad
$=3\{A(x,Jx,Jx,x)+A(x,y,y,x)+A(Jy,Jx,Jx,Jy)+A(Jy,y,y,Jy)$\smallbreak\qquad
$+2A(x,Jx,Jx,Jy)+2A(x,y,y,Jy)-2A(x,Jx,y,x)$\smallbreak\qquad
$-2A(Jy,Jx,y,Jy)-2A(x,Jx,y,Jy)-2A(x,y,Jx,Jy)\}$\smallbreak\qquad
$+3\{A(x,Jx,Jx,x)+A(x,y,y,x)+A(Jy,Jx,Jx,Jy)+A(Jy,y,y,Jy)$\smallbreak\qquad
$-2A(x,Jx,Jx,Jy)-2A(x,y,y,Jy)+2A(x,Jx,y,x)$\smallbreak\qquad
$+2A(Jy,Jx,y,Jy)-2A(x,Jx,y,Jy)-2A(x,y,Jx,Jy)\}$\smallbreak\qquad
$-\{A(x,Jx,Jx,x)+A(x,Jy,Jy,x)+A(y,Jx,Jx,y)+A(y,Jy,Jy,y)$\smallbreak\qquad
$+2A(x,Jx,Jx,y)+2A(x,Jy,Jy,y)+2A(x,Jx,Jy,x)$\smallbreak\qquad
$+2A(y,Jx,Jy,y)+2A(x,Jx,Jy,y)+2A(x,Jy,Jx,y)\}$\smallbreak\qquad
$-\{A(x,Jx,Jx,x)+A(x,Jy,Jy,x)+A(y,Jx,Jx,y)+A(y,Jy,Jy,y)$\smallbreak\qquad
$-2A(x,Jx,Jx,y)-2A(x,Jy,Jy,y)-2A(x,Jx,Jy,x)$\smallbreak\qquad
$-2A(y,Jx,Jy,y)+2A(x,Jx,Jy,y)+2A(x,Jy,Jx,y)\}$\smallbreak\qquad
$-4\{A(x,Jx,Jx,x)+A(y,Jy,Jy,y)\}$\goodbreak\medbreak\qquad
$=6A(x,y,y,x)+6A(Jx,Jy,Jy,Jx)-2A(x,Jy,Jy,x)-2A(y,Jx,Jx,y)$\smallbreak\qquad
$-12A(x,y,Jx,Jy)-12A(x,Jx,y,Jy)-4A(x,Jx,Jy,y)-4A(x,Jy,Jx,y)$\goodbreak\medbreak\qquad
$=12A(x,y,y,x)-4A(x,Jy,Jy,x)$\smallbreak\qquad
$-12A(x,y,Jx,Jy)-8A(x,Jx,y,Jy)-4A(x,Jy,Jx,y)$.\medbreak\noindent
We may now compute: \medbreak\quad
$3Q(x+Jy)+3Q(x-Jy)-Q(x+y)-Q(x-y)$ \smallbreak\qquad
$-4Q(x)-4Q(y)+4\{5\lambda(x,y)+\lambda(x,Jy)\}$ \smallbreak\quad
$=12A(x,y,y,x)-4A(x,Jy,Jy,x)$\smallbreak\qquad
$-12A(x,y,Jx,Jy)-8A(x,Jx,y,Jy)-4A(x,Jy,Jx,y)$\smallbreak\qquad
$+20(A(x,y,y,x)-A(x,y,Jy,Jx))$\smallbreak\qquad
$+4(A(x,Jy,Jy,x)+A(x,Jy,y,Jx))$\smallbreak\quad
$=32A(x,y,y,x)-8A(x,y,Jy,Jx)-8A(x,Jx,y,Jy)-8A(x,Jy,Jx,y)$.
\medbreak\noindent The desired identity now follows from the
Bianchi identity.\end{proof}

The following tensors arise naturally and play a fundamental role
in studying complex models. The tensor $A_0$ has constant
sectional curvature $+1$ and the curvature tensor of complex
projective space with the Fubini-Study metric is given by
$A_0+A_J$ where we define:
\begin{equation}\label{eqn-can-tens}
\begin{array}{l}
A_0(x,y,z,w):=\g{x}{w}\g{y}{z}-\g{x}{z}\g{y}{w}\,,\\
A_J(x,y,z,w):=\g{x}{Jw}\g{y}{Jz}-\g{x}{Jz}\g{y}{Jw}-2\g{x}{Jy}\g{z}{Jw}\,.\vphantom{\vrule
height 11pt}
\end{array}\end{equation}

The following result is due to Sato \cite{sato}; again, it was
stated in a geometrical context.

\begin{lemma}\label{sato-lemma}
Let $\mathfrak{C}=(V,\g\cdot\cdot,A,J)$ be a compatible complex
model.
\begin{enumerate}
\item If $\mathfrak{C}$ has constant holomorphic sectional
curvature $c$, then \smallbreak $A(x,y,z,w)=
\frac{c}{4}\{A_0(x,y,z,w)+A_J(x,y,z,w)\}$\smallbreak\quad
$+\frac{1}{8}\{5A(x,y,z,w)-3A(x,y,Jz,Jw)+A(x,z,Jw,Jy)$\smallbreak\quad
$-A(x,w,Jz,Jy)-A(x,Jz,w,Jy)+A(x,Jw,z,Jy)\}$. \smallbreak\item If
$\mathfrak{C}$ has constant zero holomorphic sectional curvature,
then \smallbreak $3A(x,y,z,w)+3A(x,y,Jz,Jw)$\smallbreak\noindent$
=A(x,z,Jw,Jy)-A(x,w,Jz,Jy)
-A(x,Jz,w,Jy)+A(x,Jw,z,Jy)\,.$\end{enumerate}
\end{lemma}
\begin{proof}As the holomorphic sectional curvature is constant,
 $Q(x)=c\g{x}{x}^2$. We use the identity of Lemma \ref{vanhecke-lemma} to
see: \medbreak\qquad$32A(x,y,y,x)$\smallbreak\qquad
$=3c\{\g{x}{x}+\g{y}{y}+2\g{x}{Jy}\}^2
+3c\{\g{x}{x}+\g{y}{y}-2\g{x}{Jy}\}^2$\smallbreak\qquad
$-c\{\g{x}{x}+\g{y}{y}+2\g{x}{y}\}^2-c\{\g{x}{x}+\g{y}{y}-2\g{x}{y}\}^2$\smallbreak\qquad
$-4c\g{x}{x}^2-4c\g{y}{y}^2
+4\{5\lambda(x,y)+\lambda(x,Jy)\}$\smallbreak\qquad
$=8c\{\g{x}{x}\g{y}{y}-\g{x}{y}^2+3\g{x}{Jy}^2\}
+4\{5\lambda(x,y)+\lambda(x,Jy)\}$.\medbreak\noindent We now
polarize this identity to see:
\begin{eqnarray}
&&8(A(x,y,z,w)+A(x,z,y,w))\nonumber\\
&=& 2c\{2\g{x}{w}\g{y}{z}-\g{x}{y}\g{z}{w}-\g{x}{z}\g{w}{y}\nonumber\\
&+&3\g{x}{Jy}\g{w}{Jz}+3\g{x}{Jz}\g{w}{Jy}\}\label{polarized-1}\\
&+&5\{A(x,y,z,w)+A(x,z,y,w)-A(x,y,Jz,Jw)-A(x,z,Jy,Jw)\}\nonumber\\
&+&A(x,Jy,z,Jw)+A(x,Jz,y,Jw)+A(x,Jy,Jz,w)+A(x,Jz,Jy,w)\,.\nonumber
\end{eqnarray}
Interchanging $x$ and $y$ in Equation (\ref{polarized-1}) yields:
\begin{eqnarray}
&&8(A(y,x,z,w)+A(y,z,x,w))\nonumber\\
&=&
2c\{2\g{y}{w}\g{x}{z}-\g{y}{x}\g{z}{w}-\g{y}{z}\g{w}{x}\nonumber\\
&+&3\g{y}{Jx}\g{w}{Jz}+3\g{y}{Jz}\g{w}{Jx}\}\label{polarized-2}\\
&+&5\{A(y,x,z,w)+A(y,z,x,w)-A(y,x,Jz,Jw)-A(y,z,Jx,Jw)\}\nonumber\\
&+&A(y,Jx,z,Jw)+A(y,Jz,x,Jw)+A(y,Jx,Jz,w)+A(y,Jz,Jx,w)\,.\nonumber
\end{eqnarray}
We now subtract Equation (\ref{polarized-2}) from Equation
(\ref{polarized-1}) and simplify to obtain: \medbreak\quad
$8(A(x,y,z,w)+A(x,z,y,w)-A(y,x,z,w)-A(y,z,x,w))$\smallbreak\quad$
=24 A(x,y,z,w)$\smallbreak\quad$
=2c\{3\g{x}{w}\g{y}{z}-3\g{x}{z}\g{w}{y}$\smallbreak\quad$
+3\g{x}{Jw}\g{y}{Jz}-3\g{x}{Jz}\g{y}{Jw}-6\g{x}{Jy}\g{z}{Jw}\}$\smallbreak\quad$
+10A(x,y,z,w)+5A(x,w,z,y)-5A(x,z,w,y)$\smallbreak\quad$
-10A(x,y,Jz,Jw)+A(x,Jw,Jz,y)-A(x,Jz,Jw,y)$\smallbreak\quad$
-5A(x,w,Jz,Jy)+A(x,Jy,Jz,w)-A(x,Jz,w,Jy)+A(x,Jy,Jz,w)$\smallbreak\quad$
+5A(x,z,Jw,Jy)+A(x,Jy,z,Jw)+A(x,Jw,z,Jy)+A(x,Jy,z,Jw)$\smallbreak\quad$
= 2c\{3R_0(x,y,z,w)+3R_J(x,y,z,w)\}
+15A(x,y,z,w)-9A(x,y,Jz,Jw)$\smallbreak\quad$
-3A(x,w,Jz,Jy)-3A(x,Jz,w,Jy)
+3A(x,Jw,z,Jy)+3A(x,z,Jw,Jy)$.\medbreak\noindent The desired
results now follow.\end{proof}

We begin the proper study of the complex Jacobi operator by
examining the condition $\ja(\cdot)=0$.

\begin{lemma}\label{lem-2.3}
Let $\mathfrak{C}=(V,\g\cdot\cdot,J,A)$ be a complex model.
\begin{enumerate}\item The following
conditions are equivalent:
\begin{enumerate}
\item $\ja(\pi)=0$ for all $\pi\in\mathbb{CP}(V,J)$. \item
${\mathcal{A}}(x,y)=-{\mathcal{A}}(Jx,Jy)$ for all $x,y$. \item
${\mathcal{A}}(\pi)=0$ for all $\pi\in\mathbb{CP}(V,J)$. \item
$\mathcal{A}(Jx,y)z={\mathcal{A}}(x,Jy)z={\mathcal{A}}(x,y)Jz$ for
all $x,y,z$.
\end{enumerate}
\item If any of the conditions in (1) are satisfied, then:
\begin{enumerate}
\item $\mathfrak{C}$ is compatible. \item $\mathfrak{C}$ is Ricci
flat and $\star$-Ricci flat.
\end{enumerate}\end{enumerate}
\end{lemma}
\begin{proof} Suppose Condition (1a) holds. Then $\mathfrak{C}$ is compatible. Furthermore,
$$Q(x)=R(x,Jx,Jx,x)=R(x,Jx,Jx,x)+R(x,x,x,x)=\langle\mathcal{J}(\pi_x)x,x\rangle=0\,,$$
so $\mathfrak{C}$ has constant holomorphic sectional curvature 0.
Thus Lemma \ref{sato-lemma} (2) applies and we show that Condition
(1b) holds by computing:
\begin{eqnarray*}
&&3A(x,y,z,w)+3A(x,y,Jz,Jw)\\
&=&A(x,z,Jw,Jy)-A(x,w,Jz,Jy)-A(x,Jz,w,Jy)+A(x,Jw,z,Jy)\\
&=&\langle\{\ja(\pi_{z+Jw})-\ja(\pi_z)-\ja(\pi_{Jw})\}x,Jy\rangle=0\,.
\end{eqnarray*}

Suppose Condition (1b) holds. We establish Condition (1c) by
computing:
$$\mathcal{A}(\pi_x)=\mathcal{A}(x,Jx)=-\mathcal{A}(Jx,JJx)=\mathcal{A}(Jx,x)=-\mathcal{A}(x,Jx)
=-\mathcal{A}(\pi_x)\,.$$

Suppose Condition (1c) holds. Then $\mathfrak{C}$ is compatible.
Again Lemma \ref{sato-lemma} is applicable. We set $w=Jz$ in Lemma
\ref{sato-lemma} (2) to show Condition (1a) holds by computing:
\begin{eqnarray*}
0&=&6\langle\mathcal{A}(\pi_z)x,y\rangle=6A(x,y,z,Jz)\\
&=&-2A(x,z,z,Jy)-2A(x,Jz,Jz,Jy)=-2\langle\mathcal{J}(\pi_z)x,Jy\rangle\,.
\end{eqnarray*}

We have shown that Conditions (1a), (1b), and (1c) are equivalent.
Suppose that Condition (1d) holds. We show that Condition (1c)
holds by computing:
$$
\mathcal{A}(\pi_x)=\mathcal{A}(x,Jx)
=\mathcal{A}(Jx,x)=-\mathcal{A}(x,Jx)=-\mathcal{A}(\pi_x)\,.
$$

Finally suppose that Condition (1a) holds. We must show Condition
(1d) holds. Since Condition (1a) implies Condition (1b), we have
${\mathcal{A}}(x,y)=-{\mathcal{A}}(Jx,Jy)$. Thus
$$A(Jx,y)z=-A(JJx,Jy)z={\mathcal{A}}(x,Jy)z\,.$$
Since $\ja(\pi_y)=0$ and since $\mathfrak{C}$ is compatible,
$$0=A(x,y,y,w)+A(x,Jy,Jy,w)=
A(x,y,y,w)+A(Jx,y,y,Jw)\,.
$$
Polarize this identity to see
\begin{equation}\label{1-equation}
0=A(x,y,z,w)+A(x,z,y,w)+A(Jx,y,z,Jw)+A(Jx,z,y,Jw)\,.
\end{equation}
Since ${\mathcal{A}}(x,w)=-{\mathcal{A}}(Jx,Jw)$,
$A(Jx,Jw,y,z)=-A(x,w,y,z)$. We may therefore use the First Bianchi
identity to see:
\begin{eqnarray}
0&=&A(Jx,y,z,Jw)+A(Jx,Jw,y,z)+A(Jx,z,Jw,y)\nonumber\\
&+&A(x,y,z,w)+A(x,w,y,z)+A(x,z,w,y)\label{2-equation}\\
&=&A(Jx,y,z,Jw)-A(Jx,z,y,Jw)+A(x,y,z,w)-A(x,z,y,w)\,.\nonumber
\end{eqnarray}
Adding (\ref{1-equation}) and (\ref{2-equation}) we get
$$0=2A(x,y,z,w)+2A(Jx,y,z,Jw)\,.$$
Replacing $w$ by $Jw$ and changing the order then shows Condition
(1a) implies Condition (1d) since:
$$A(x,y,Jw,z)=A(Jx,y,w,z)\,.$$

This completes the proof of Assertion (1). It is clear that (1a)
implies $\mathfrak{C}$ is compatible, as we said before. Assume
the conditions of Assertion (1) hold. We may then compute:
\begin{eqnarray*}
&&2\rho(x,y)=\sum_{i=1}^m\{A(e_i,x,y,e_i)+A(Je_i,x,y,Je_i)\}
=\sum_{i=1}^m\langle\mathcal{J}(\pi_{e_i})x,y\rangle)=0,\\
&&\rho^\star(x,y)=\sum_{i=1}^mA(x,Je_i,Jy,e_i)=\sum_{i=1}^mA(x,e_i,JJy,e_i)=\rho(x,y)=0\,.
\end{eqnarray*}
where $\{e_1,\dots,e_m\}$ forms an orthonormal basis. This
completes the proof of the Lemma.\end{proof}

\begin{proof} We now establish Theorem \ref{thm-1.5}.
Since the dimension of $V$ is a multiple of $4$, we may choose
almost complex structures $J,K$ such that $JK+KJ=0$. Consider the
algebraic curvature tensor $A:=A_K-A_{JK}$. Note that the Jacobi
operator is given by
\[
\ja(x)y=3\g{y}{Kx}Kx-3\g{y}{JKx}JKx\,,
\]
and, since $\ja(x)Kx=Kx$ for any unit vector $x$, $A\neq 0$.
However, the complex Jacobi operator vanishes identically:
\[
\begin{array}{rcl}
\ja(\pi_x)&=&3\g{y}{Kx}Kx-3\g{y}{JKx}JKx\\
\noalign{\smallskip} &&\!\!\!\!+3\g{y}{KJx}KJx-3\g{y}{JKJx}JKJx\\
\noalign{\smallskip} &=&0\,.
\end{array}
\]
We may now apply Lemma \ref{lem-2.3} to see that
$\mathcal{A}(\pi)$ also vanishes identically as well.\end{proof}

We conclude this section by putting things in a slightly different
invariant framework. Let $\mathfrak{A}(V)$ be the vector space of
all algebraic curvature tensors on $V$. We use
$\langle\cdot,\cdot\rangle$ to define a natural inner product on
$\mathfrak{A}(V)$ by setting:
\[
\langle A_1,A_2\rangle:=\sum_{i,j,k,l}
A_1(e_i,e_j,e_k,e_l)A_2(e_i,e_j,e_k,e_l)\,;
\]
this is independent of the particular orthonormal basis $\{e_i\}$
chosen. Consider the following subspaces \cite{gray}:
\[
\begin{array}{rcl}
\mathfrak{A}_1(V,J)&=&\{A\in \mathfrak{A}(V):
A(x,y,z,w)=A(Jx,Jy,z,w)\}\,,\\
\noalign{\smallskip}
\mathfrak{A}_2(V,J)&=&\{A\in \mathfrak{A}(V):A(x,y,z,w)=A(Jx,Jy,z,w)\\
\noalign{\smallskip}    && \qquad\qquad\qquad\qquad +A(Jx,y,Jz,w)+A(Jx,y,z,Jw)\},\\
\noalign{\smallskip} \mathfrak{A}_3(V,J)&=&\{A\in \mathfrak{A}(V):
A(x,y,z,w)=A(Jx,Jy,Jz,Jw)\}\,.
\end{array}
\]
Note that $\mathfrak{A}_1(V,J)$ is the space of algebraic
curvature tensors which verify the K{\"a}hler identity and
$\mathfrak{A}_3(V,J)$ is the space of compatible curvature
tensors. We have
$$\mathfrak{A}_1(V,J)\subset \mathfrak{A}_2(V,J)\subset
\mathfrak{A}_3(V,J)\,.$$ Denote by
$\mathfrak{A}_2^\bot(V,\langle\cdot,\cdot\rangle,J)$ the
orthogonal complement of $\mathfrak{A}_2(V,J)$ in
$\mathfrak{A}_3(V,J)$.

\begin{lemma}\label{lem-2.4}
Let $\mathfrak{C}=(V,\langle\cdot,\cdot\rangle,A)$ be a complex
model. The following assertions are equivalent:
\begin{enumerate}
\item $\mathcal{J}(\pi)=0$ for all $\pi\in\mathbb{CP}(V,J)$. \item
$A\in\mathfrak{A}_2^\bot(V,\langle\cdot,\cdot\rangle,J)$.
\end{enumerate}
\end{lemma}

\begin{proof} Projection in $\mathfrak{A}_2(V,J)$ restricted to $\mathfrak{A}_3(V,J)$
is given by
\[
\begin{array}{rcl}
\mathcal{P}_2(A)(x,y,z,w)&=&\frac{1}{2}\{A(x,y,z,w)+ A(Jx,Jy,z,w)\\
 \noalign{\smallskip}&&\quad+A(Jx,y,Jz,w)+A(Jx,y,z,Jw)\}\,.
\end{array}
\]
$\mathcal{P}_2$ is an involutive isometry and
$\mathfrak{A}_2^\bot(V,\langle\cdot,\cdot\rangle,J)$ is the
$(-1)$-eigenspace of $\mathcal{P}_2$ by \cite{tricerri-vanhecke}. For $A\in\mathfrak{A}_3(V,J)$
compute
\[
\mathcal{P}_2(A)(x,y,z,w)+\mathcal{P}_2(A)(Jx,Jy,z,w)=A(x,y,z,w)+
A(Jx,Jy,z,w)\,.
\]
From here, it is easy to verify that
\[
\mathfrak{A}_2^\bot(V,\langle\cdot,\cdot\rangle,J)=\{A\in\mathfrak{A}(V):
A(x,y,z,w)=-A(Jx,Jy,z,w)\}\,.
\]
The desired result now follows from Lemma \ref{lem-2.3}.
\end{proof}

\section{Geometrical Results}\label{sect-geo}

We begin our study of the geometrical context by recalling several
well known results. We refer to Gray \cite{gray} for the proof of
Assertions (1) and (2), see also \cite{gray-1969}, and to Yano
\cite{yano} for the proof of Assertion (3) in the following Lemma:

\begin{lemma}\label{lem-3.1}
Let $\mathcal{C}=(M,g,J)$ be an Hermitian manifold and let
$\mathfrak{C}=\mathfrak{C}(\mathcal{C},P)$ be the almost complex
model determined by $\mathcal{C}$ at a point $P\in M$. Then:
\begin{enumerate}
\item If $\mathcal{C}$ is Hermitian, then
\begin{eqnarray*}
&&R(x,y,z,w)+R(Jx,Jy,Jz,Jw)= R(Jx,Jy,z,w)+R(x,y,Jz,Jw)\\
&&\qquad+R(Jx,y,Jz,w)+R(x,Jy,z,Jw)+R(Jx,y,z,Jw)+R(x,Jy,Jz,w)\,.
\end{eqnarray*}
\item If $\mathcal{C}$ is nearly K\"ahler, then $\mathcal{C}$ is
compatible and
\begin{eqnarray*}
&&R(x,y,z,w)+R(Jx,Jy,Jz,Jw)= R(Jx,Jy,z,w)+R(x,y,Jz,Jw)\\
&&\qquad+R(Jx,y,Jz,w)+R(x,Jy,z,Jw)+R(Jx,y,z,Jw)+R(x,Jy,Jz,w)\,.
\end{eqnarray*}
\item If $\mathcal{C}=(M,g,J)$ is an almost K\"ahler manifold,
then $\|\nabla J\|^2=2(\tau^\star -\tau)$.
\end{enumerate}\end{lemma}

Let $\mathcal{C}_i$ be almost Hermitian manifolds. We suppose
given a complex isometry
$\theta:(T_PM_1,g_1,J_1)\rightarrow(T_QM_2,g_2,J_2)$. We let
$V=T_PM_1$, $\langle\cdot,\cdot\rangle=g_1$, $J=J_1$, and
$A:=R_1-\theta^*R_2$. Let
$$\mathfrak{C}=\mathfrak{C}(\mathcal{C}_1,P,\mathcal{C}_2,Q,\theta):=(V,\langle\cdot,\cdot\rangle,J,A)\,.$$
Theorem \ref{complex-dif-theorem} will follow from the following
Lemma:

\begin{lemma}\label{lem-3.2} Adopt the notation established above. Assume that $\mathcal{C}_i$ are Hermitian
or nearly K{\"a}hler manifolds. The following assertions are
equivalent:
\begin{enumerate}
\item $\mathcal{J}(\pi)=0$ for all $\pi\in\mathbb{CP}(V,J)$. \item
$\mathcal{R}(\pi)=0$ for all $\pi\in\mathbb{CP}(V,J)$. \item
$A=0$.
\end{enumerate}
\end{lemma}

\begin{proof}
Assume that either Condition (1) or Condition (2) holds; these are
equivalent by Lemma \ref{lem-2.3}. Since the curvature tensors
$A_{\mathcal{C}_i}$ satisfy the identity of Lemma \ref{lem-3.1},
so does their difference. We use the relations provided by Lemma
\ref{lem-2.3} to show that $A=0$ by computing:
\begin{eqnarray*}
0&=&R(x,y,z,w)+R(Jx,Jy,Jz,Jw)- R(Jx,Jy,z,w)-R(x,y,Jz,Jw)\\
&-&R(Jx,y,Jz,w)-R(x,Jy,z,Jw)-R(Jx,y,z,Jw)-R(x,Jy,Jz,w)\\
&=&R(x,y,z,w)+R(x,y,z,w)- R(JJx,y,z,w)-R(JJx,y,z,w)\\
&-&R(JJx,y,z,w)-R(JJx,y,z,w)-R(JJx,y,z,w)-R(JJx,y,z,w)\\
&=&8R(x,y,z,w)\,.
\end{eqnarray*}
Conversely, of course, if $A=0$, then
$\mathcal{J}(\pi)=\mathcal{R}(\pi)=0$ for all
$\pi\in\mathbb{CP}(V,J)$.
\end{proof}

\begin{proof} We now prove Theorem \ref{thm-1.3}.
We use Lemma \ref{lem-2.3} to see that $\mathcal{M}$ is both Ricci
flat and $\star$-Ricci flat. Hence $\tau=\tau^\star=0$. Therefore,
by Lemma \ref{lem-3.1} (3), $\nabla J=0$ and the manifold is
K{\"a}hler. This implies the almost complex structure is in fact
integrable so $\mathcal{C}$ is Hermitian. The desired conclusion
now follows from Theorem \ref{complex-dif-theorem}.\end{proof}

\begin{proof} We now prove Theorem \ref{thm-1.6}. Our construction is motivated by the construction of Theorem
\ref{thm-1.5} and is based on work of Sato \cite{sato-2003}. Let
$m=4n$. Let $(\mathbb{CP}^{2n},g,J)$ be complex projective space
with the Fubini-Study metric $g$ and usual complex structure
$J_0$; this is a K\"ahler manifold. The canonical embedding of
$\mathbb{C}^{2n}\subset\mathbb{C}^{2n+1}$ defines an isometric
embedding of $\mathbb{CP}^{2n-1}$ in $\mathbb{CP}^{2n}$. Let
$M:=\mathbb{CP}^{2n}-\mathbb{CP}^{2n-1}$. Let $\mathcal{H}$ be the
fiber bundle of all unitary quaternion structures
$\{J_1,J_2,J_3\}$ on the tangent bundle of $M$ which satisfy
$J_1=J$. Since $M$ is contractable, $\mathcal{H}$ is a trivial
fiber bundle so we can define a global quaternion structure
$\{J_1,J_2,J_3\}$ on $TM$ so that $J=J_1$. This is, of course,
just the usual twistor construction.

Let $\mathcal{C}:=(M,g,J_2)$. Let
$$\Theta:x\rightarrow(1+J_2)/\sqrt{2}\,.$$
This defines an isometry of $T_PM$ with $\Theta J_2=J_2\Theta$.
Furthermore
$$\Theta J_1=-J_3\Theta\quad\text{and}\quad\Theta J_3=J_1\Theta\,.$$
The curvature tensor of the Fubini-Study metric is given by
$R_0+R_{J_1}$. Let $x$ be a unit tangent vector. We use the
defining relations of Equation (\ref{eqn-can-tens}) to see that:
$$\mathcal{J}_R(x)y=\left\{\begin{array}{rll}
0&\text{if}&y\in\operatorname{Span}\{x\},\\
4y&\text{if}&y\in\operatorname{Span}\{J_1x\},\\
y&\text{if}&y\perp\operatorname{Span}\{x,J_1x\}\,.
\end{array}\right.$$
As $\Theta^*R=R_0+R_{J_3}$ and as $J_1x\perp J_3x$, $\Theta^*R\ne
R$. Since
$\mathcal{J}_R(\pi_x)=\mathcal{J}_R(x)+\mathcal{J}_R(J_2x)$,
$$\mathcal{J}_R(\pi_x)y=\left\{\begin{array}{rll}
4y&\text{if}&y\in\operatorname{Span}\{x,J_2x\},\\
5y&\text{if}&y\in\operatorname{Span}\{J_1x,J_1J_2x=J_3x\},\\
2y&\text{if}&y\perp\operatorname{Span}\{x,J_1x,J_2x,J_3x\}\,.
\end{array}\right.$$
Since $J_1$ and $J_3$ play symmetric roles in this identity,
$\mathcal{J}_{\Theta^*R}(\pi_x)=\mathcal{J}_{R}(\pi_x)$ as
desired. Lemma \ref{lem-2.3} now shows
$\mathcal{R}_{\Theta^*R}(\pi_x)=\mathcal{R}_R(\pi_x)$ as well.
\end{proof}

\section{Conformal and almost K\"ahler geometry}\label{sect-add-res}

\begin{theorem}
Let $\mathcal{C}=(M,g,J)$ be an almost K{\"a}hler manifold such that
there exists a K{\"a}hler metric $\tilde g$ on $(M,J)$ such that
${\mathcal{R}}(\pi)=\tilde{\mathcal{R}}(\pi)$ for all
$\pi\in\mathbb{CP}(TM,J)$ (equivalently,
$\ja_R(\pi)=\mathcal{J}_{\tilde R}(\pi)$ for all
$\pi\in\mathbb{CP}(TM,J)$). Then $\mathcal{C}$ is K{\"a}hler and
$R=\tilde R$.
\end{theorem}
\begin{proof} Suppose ${\mathcal{R}}(x,Jx)=\tilde{\mathcal{R}}(x,Jx)$. Set $\bar R:=R-\tilde R$. Then
$\bar{\mathcal R}(\pi)=0$ for all $\pi$. Since $\tau^\star=\tau$
for any K{\"a}hler manifold and $\bar{R}$ is Ricci flat and
$\star$-Ricci flat by Lemma \ref{lem-2.3} (2), one has $\|\nabla
J\|^2=0$ by Lemma \ref{lem-3.1} (3); hence $\nabla J=0$ and the
original manifold $\mathcal{C}$ is indeed K{\"a}hler. That the
curvature tensors are equal follows from
Theorem~\ref{complex-dif-theorem}.\end{proof}

\begin{theorem}
Let $\mathcal{C}:=(M,g,J)$ and $\mathcal{C}^\alpha:=(M,e^\alpha
g,J)$ be conformally equivalent almost Hermitian manifolds. If
$\mathcal{J}_{\mathcal
C}(\pi_x)=\mathcal{J}_{\mathcal{C}^\alpha}(\pi)$ for all
$\pi\in\mathbb{CP}(TM,J)$, then $R=R^\alpha$.
\end{theorem}
\proof Let $\sigma(g,J)$ be projection on
$\mathfrak{A}_2^\bot(g,J)$. A priori this projection depends on
the choice of $g$. It is, however conformally invariant, i.e.
$\mathfrak{A}_2^\bot(g,J)=\mathfrak{A}_2^\bot(e^\alpha g,J)$ and
$\sigma(g,J)=\sigma(e^\alpha g,J)$. Set $\sigma:=\sigma(g,J)$.
Furthermore, see \cite{tricerri-vanhecke},
$\sigma(R)=\sigma(R^\alpha)$; thus $\sigma(R-R^\alpha)=0$. Since
by hypothesis $(\mathcal{J}_{\mathcal
C}-\mathcal{J}_{\mathcal{C}^\alpha})(\pi_x)=0$ for all $x$, we
have as desired that $R=R^\alpha$, using
Theorem~\ref{complex-dif-theorem}. \qed

\section*{Acknowledgments}
The research of M. Brozos-V\'azquez was partially supported by
project MTM2005-08757-C04-01 (Spain) and a FPU grant. The research
of E. Garc{\'\i}a-R{\'\i}o was partially supported by project BFM2003-02949
(Spain). The research of both M. Brozos-V\'azquez and P. Gilkey
was partially supported by the Max Planck Institute for
Mathematics in the Sciences (Leipzig, Germany).

\end{document}